\documentclass[draft,12pt]{amsart}
\usepackage{latexsym}
\usepackage{amsfonts}
\usepackage{amsmath}
\usepackage{fancyheadings}
\usepackage{amssymb}
\usepackage{graphics}
\usepackage[OT2,OT1]{fontenc}
\usepackage{wrapfig}
\usepackage{color}

\newtheorem{theorem}{Theorem}[section]

\newtheorem{lemma}[theorem]{Lemma}
\newtheorem{corollary}[theorem]{Corollary}

\theoremstyle{definition}
\newtheorem{example}[theorem]{Example}

\newcommand\Z{\mathbb Z}

\newcommand\F{\mathbb F}

\newcommand{\der}[2]{#1^{(#2)}}

\newcommand{\lcs}[2]{\gamma_{#2}(#1)}

\newcommand{\sym}{\sf S}

\newcommand{\gen}[1]{\bigl<#1\bigr>}

\newcommand{\ord}[1]{|\,#1\,|}

\newcommand\rc[3]{[#1,\,{}^{#3}\!#2]}
\newcommand\e{\equiv}
\renewcommand{\d}{d}

\renewcommand{\leq}{\leqslant}
\renewcommand{\geq}{\geqslant}

\begin{document}

\title{The derived series of a finite $p$-group}
\author{Csaba Schneider}
\address{Informatics Laboratory\\ 
Computer and Automation
Research Institute\\
The Hungarian Academy of Sciences\\
1518 Budapest Pf.\ 63\\
Hungary}
\email{csaba.schneider@sztaki.hu\protect{\newline} {\it WWW:}
www.sztaki.hu/$\sim$schneider}

\begin{abstract}
Let $G$ be a finite $p$-group, and let 
$\der Gd$ denote the $d$-th term of the derived series of $G$. 
We show, for $p\geq 5$, 
that $\der Gd\neq 1$ implies 
$\log_p\ord G\geq 2^d+3d-6$, and hence we improve a recent result by Mann.
\end{abstract}

\date{\today}
\keywords{finite $p$-groups, derived subgroups, derived quotients, derived
  series, associated Lie ring}
\subjclass[2000]{20D15, 20-04}

\maketitle
\section{Introduction}

In a group $G$, let $\der Gd$ denote the $d$-th term of the derived series (as
usual, we write $\der G0=G$, $\der G1=G'$, $\der G2=G''$).
The main goal of this article is to give a new lower bound for the order of a
$p$-group in which $\der Gd\neq 1$.

\begin{theorem}\label{main}
Let  $p$ be a prime such that $p\geq 5$, 
and let $G$ be a 
finite $p$-group. If $\der Gd\neq 1$ then $\log_p\ord G\geq 2^d+3d-6$. 
\end{theorem}

The problem of bounding the order of a finite $p$-group in terms of its
derived length has a long 
history starting from Burnside's papers~\cite{burn1,burn2}. 
It is well-known that non-abelian $p$-groups have order at least 
$p^3$ and non-metabelian $p$-groups have order 
at least $p^6$ (for $p=2,\ 3$, the sharp bound is $p^7$).
In his 1934 paper~\cite{Hall34}, 
P.~Hall noticed, for $i\geq 0$, that $\der G{i+1}\neq 1$ implies 
$\log_p|\der Gi/\der G{i+1}|\geq 2^i+1$. 
He used  this lower bound  to prove that
\begin{equation}\label{hallth}
\der Gd\neq 1\quad\mbox{implies}\quad \log_p|G|\geq 2^d+d.
\end{equation} 

Recently, there has been some interest in obtaining  sharper versions
of Hall's theorem.
Mann~\cite{Mann} showed that if $G$ is a $p$-group then $\der Gd\neq 1$ 
implies $\log_p|G|\geq 2^d+2d-2$, and Theorem~\ref{main} in this paper improves
the linear term of Mann's bound. Another direction in this research is to
consider the assertion~\eqref{hallth} for a fixed $d$.
Hall's bound
is, in general, sharp for $d=1,\ 2$ (we mentioned above that, for $p=2,\ 3$
and  $d=2$, the sharp bound is $p^7$), but it is not
sharp already 
for $d=3$. In this case, Hall's lower bound is 
$p^{11}$, but, for $p\geq 3$, the smallest 
$p$-groups in which $\der G3\neq 1$ have
order $p^{14}$ 
(see~\cite{Susan,se-rmfncs,susan-mfn}) 
and the smallest such 2-groups
have order $2^{13}$~\cite{susan-mfn}. 
For $p\geq 5$ and $d\geq 3$, the smallest known 
groups with $\der Gd\neq 1$ have
composition length $2^{d+1}-2$ (see~\cite{se-rmfncs}), and so 
there is a gap to be filled.

Similar questions arise in the context of soluble groups.
Glasby~\cite{gl1,gl2} 
studied soluble groups with a fixed derived length and smallest possible 
composition length.

We prove Theorem~\ref{main} by studying the abelian and metabelian quotients of
the derived series of a finite $p$-group and obtain several results that are
interesting in their own right. By Hall's lemma referred to above, we have, in
a $p$-group $G$, that $\der G{d+1}\neq 1$ implies $\log_p|\der Gd/\der
G{d+1}|\geq 2^d+1$. In a finite $p$-group $G$ we call $\der Gd/\der G{d+1}$ a
{\em small derived quotient} if 
\begin{equation}\label{hallsharp}
\der G{d+1}\neq 1\quad\mbox{and}\quad \log_p|\der Gd/\der G{d+1}|=2^d+1.
\end{equation}
We do not 
know whether, for an arbitrary $p$ and $d$, the quotient
$\der Gd/\der G{d+1}$ can be small in a $p$-group $G$.
Groups  in which $G/G'$ or 
$G'/G''$ 
is small
are easy to construct; odd-order groups with small second derived 
quotients $G'/G''$ were 
characterised in~\cite{sch03}. The Sylow
2-subgroup $P$ of the symmetric group $\sym_{2^\d}$ of rank $2^\d$ satisfies
$\log_2\ord{\der{P}{\d-2}/\der{P}{\d-1}}=2^{\d-2}+1$ and $\der{P}{\d-1}\neq 1$
(see \cite[Lemma~(II.7)]{lg}). Therefore small derived quotients do indeed
exist. 
On the other hand, for $p\geq 3$ and $d\geq 2$, 
we do not have examples in which $\der Gd/\der G{d+1}$ is small.
We do, however, know that the structure of a finite $p$-group 
in which $\der Gd/\der G{d+1}$ is small
must be
severely restricted. Blackburn in~\cite{black87} 
attributes a theorem to P.\ Hall that, for odd primes, if $G'/G''$ is small, 
then $|G''|=p$; see~\cite{black87}. In the same article,
Blackburn himself shows that the same condition for 2-groups implies that
$G''$ is a two-generator abelian group and is nearly homocyclic.
Using Bokut's ideas~\cite{bokut71}, this result was generalised
in my PhD Thesis~\cite{csthesis}, where I showed that 
if $p\neq 3$ and~\eqref{hallsharp} holds with $d\geq 1$, then $\der G{d}$ has 
nilpotency class at most~3.
Further, Mann~\cite{Mann} showed that a finite $p$-group can have at most two
small derived quotients.

It is explained in Section~\ref{lie} that if~\eqref{hallsharp} holds in a finite
$p$-group $G$ then $\der Gd/[\der Gd,G]$ must be elementary abelian with 
rank at most~2. In this article, I prove the following theorem, which claims
that $p^2$ is usually not possible.

\begin{theorem}\label{th1}
Suppose that $p$ is an odd prime, $d\geq 1$, and $G$ is a finite $p$-group
such that $\der G{d+1}\neq 1$. If $\log_p|\der Gd/\der G{d+1}|=2^d+1$,
then $|\der Gd/[\der Gd,G]|=p$. 
\end{theorem}

Theorem~\ref{th1} is a generalisation of~\cite[Theorem~1.1]{sch03}. Its 
proof will be presented in Section~\ref{3.4} using the Lie
algebra constructed in Section~\ref{liealg}.
In a separate article~\cite{sch2}, I use Theorem~\ref{th1} to show, 
for odd $p$, that
two distinct small derived quotients can occur only in a group with order
$p^6$. This gives an improvement of Mann's related result in~\cite{Mann}.

Theorem~\ref{th1} is used to obtain the constant term ``6'' in
Theorem~\ref{main}. In order to obtain the coefficient ``3'' of the linear
term, we study metabelian quotients of the derived series, that is, we 
consider the factor groups
$\der Gi/\der G{i+2}$. This way we are able to prove a better lower bound for the order of $G$ than those obtained by
considering only the abelian quotients $\der Gi/\der G{i+1}$. 

The structure of this paper is as follows. In Section~\ref{tools} we summarise the
technical tools that we need in the proofs of the main theorems. We prove
Theorem~\ref{th1} in Sections~\ref{lie}--\ref{3.4}. Many of the 
proofs 
in this paper use the Lie ring method. In order to verify Theorem~\ref{th1} we
define a filtration on $\der Gd$ and consider the graded Lie ring
associated with this filtration. We identify a small subring in this Lie ring,
and relatively simple, though rather lengthy, computation yields
Theorem~\ref{th1}. Then, in Section~\ref{secth2}, we study the metabelian
quotients of the derived series of a finite $p$-group, and prove
the somewhat technical Theorem~\ref{p3}. Again, the hard part is to identify
a small subring in the Lie ring associated with the lower central series of
$G$. Once this subring is found, elementary, but somewhat cumbersome,
calculations leads to Theorem~\ref{p3}. Finally, in Section~\ref{last} we
present a proof for Theorem~\ref{main}.

\section{A toolbox}\label{tools}

In this section we review a couple of well-known results that are used in this
paper. 
The construction of the Lie ring in Section~\ref{liealg} is based on some
simple properties of group commutators. 
In addition to the best-known commutator identities (see for instance~\cite[Proposition~1.1.6]{l-g}), the following 
collection formula will often be used: if $x$ and $y$ are elements of a finite
$p$-group then
\begin{equation}\label{coll}
[x^p,y]\e[x,y]^p\bmod {H'}^p\lcs{H}{p}
\quad\mbox{where}\quad
H=\gen{x,[x,y]};
\end{equation}
see~\cite[Proposition~1.1.32]{l-g}. 
The Hall-Witt identity will occur in a lesser-known form: if $x$, $y$, and $z$
are elements of a group then
\begin{equation}\label{hw}
[x,y,z^x][z,x,y^z][y,z,x^y]=1;
\end{equation}
see~\cite[Proposition~1.1.6(v)]{l-g}.

The following result was proved by P.\ Hall~\cite{Hall34}.

\begin{lemma}\label{halllem}
Suppose that $H$ is a non-abelian normal subgroup in a finite $p$-group $G$
such that $H\leq \lcs Gi$. Then $|H/H'|\geq p^{i+1}$ and $|H|\geq p^{i+2}$. 
\end{lemma}

For the proof of the following lemma, see~\cite[Lemma~2.1]{blackburn54}.

\begin{lemma}\label{cyclemma}
If $G$ is a group and $H$ is a normal subgroup with cyclic quotient,
then $G'= [G,H]$. 
\end{lemma}

If $U$ and $V$ are subgroups in a group and $n$ is a natural
number then let $\rc{U}{V}{n}$ denote the left-normed commutator subgroup
$$
\rc{U}{V}{n}=[U,\underbrace{V,\ldots,V}_{n\rm\ copies}].
$$
The following result is easy to verify by induction on $i$.

\begin{lemma}\label{leftnormed}
If $A,\ B$ are normal subgroups of a group $G$, then
$[A,\lcs Bi]\leq \rc ABi$. 
\end{lemma}

\section{A filtration on $\der Gd$}\label{lie}

In this section we start the preparation for the proof of Theorem~\ref{th1}. 
Let $G$ be a finite $p$-group. 
Let us define a filtration on $\der Gd$ as follows:
$$
N_{2^\d}=\der Gd\quad\mbox{and}\quad N_{k+1}=[N_k,G]\quad\mbox{for}\quad k\geq 2
^\d.
$$
Using the notation introduced in the previous section, 
we may write 
$N_k=\rc{\der Gd}{G}{k-2^\d}$. The following lemma shows, in the
terminology of~\cite{Shalev95}, that
the $N_i$ form 
a strongly central filtration. 

\begin{lemma}\label{strcent}
The filtration $\{N_k\}_{k\geq 2^\d}$ satisfies the property
$[N_j,N_k]\leq N_{j+k}$ for each $j,\ k\geq 2^\d$. 
\end{lemma}
\begin{proof}
As $N_k$ is contained in $\lcs Gk$, and, 
for a subgroup $H$, the commutator subgroup $[H,\lcs Gi]$ 
is contained in $\rc HG{i}$ (see Lemma~\ref{leftnormed}), 
we obtain
$$[N_j,N_k]
\leq [N_j,\lcs Gk]\leq \rc{N_j}{G}k=N_{j+k}.
$$
\end{proof}

\begin{lemma}\label{2poss}
Suppose that the derived quotient $\der Gd/\der G{d+1}$ is small in a finite
$p$-group $G$. Then the following is true.
\begin{enumerate}
\item[(a)] Exactly one of the following must hold:
\begin{enumerate}
\item[(i)] $|N_i/N_{i+1}|=p$ for
  $i\in\{2^d,\ldots,2^{d+1}\}$, and $\der G{d+1}=N_{2^{d+1}+1}$;
\item[(ii)] $|\der Gd/N_{2^d+1}|=p^2$, $|N_i/N_{i+1}|=p$ for
  $i\in\{2^d+1,\ldots,2^{d+1}-1\}$, and $\der G{d+1}=N_{2^{d+1}}$.
\end{enumerate}
\item[(b)] $N_i/N_{i+1}$ is elementary abelian for $i\geq 2^d$.
\item[(c)] $N_{2^d+1}=\der Gd\cap \lcs G{2^d+1}$.
\end{enumerate}
\end{lemma}
\begin{proof}
(a) As $\der Gd\leq\lcs G{2^d}$, we obtain
$$
\der G{d+1}
=[\der Gd,\der Gd]\leq [
\der Gd,\lcs{G}{2^\d}]\leq \rc{\der Gd}{G}{2^\d}=N_{2^{d+1}},
$$
therefore  we have
the following chain of $G$-normal subgroups:
\begin{equation}\label{chain}
\der Gd>[\der Gd,G]>[\der Gd,G,G]>\cdots>\rc {\der Gd}{G}{2^\d}\geq \der G{d+1}.
\end{equation}
Counting number of non-trivial factors of this chain, we obtain that $\der Gd/[\der Gd,G]$ has order at most $p^2$.
If $\der Gd/[\der Gd,G]$ is cyclic, then, by Lemma~\ref{cyclemma},
the subgroup $\der G{d+1}$ coincides with $[\der Gd,[\der Gd,G]]$, and so 
$$
\der G{d+1}
=[\der Gd,\der Gd]=[\der Gd,[\der Gd,G]]\leq [
\der Gd,\lcs{G}{2^\d+1}]\leq \rc{\der Gd}{G}{2^\d+1}=N_{2^{d+1}+1}.
$$
Thus, in this case, we
obtain the following modified chain:
\begin{equation}\label{newchain}
\der Gd>[\der Gd,G]>[\der Gd,G,G]>\cdots>\rc{\der Gd}{G}{2^\d}>\rc{\der Gd}{G}{2
^d+1}\geq\der G{d+1}.
\end{equation}
If the first quotient in these chains has order $p$, then this quotient is
cyclic, and so~\eqref{newchain} must hold. 
In this case, counting the non-trivial factors in~\eqref{newchain}, we find 
that case~(i) must be valid. 
Now suppose that the 
first quotient $\der Gd/[\der Gd,G]$ has order $p^2$. 
Then chain~\eqref{newchain} is too long, and
so $\der Gd/[\der Gd,G]$ must be elementary abelian. 
As before, we count the number of factors
in~\eqref{chain} and find that case~(ii) must hold.

(b) In the previous paragraph we obtained that $\der Gd/[\der Gd,G]$ is
elementary abelian, and  we use
induction to show that so are the other factors $N_i/N_{i+1}$.
Let
us suppose that $N_{k-1}/N_{k}$ is elementary abelian 
for some $k\geq 2^\d+1$. 
Consider
a generator $a=[x,y]$ of $N_k$ where $x\in N_{k-1}$ and $y\in G$. Then
$[x,y]^p\e[x^p,y]$ modulo $N_{k+1}$ (equation~(\ref{coll})) 
and $x^p$ lies in $N_k$ by the induction
hypothesis. Hence we obtain that $[x^p,y]\in N_{k+1}$, and so
$[x,y]^p\in N_{k+1}$ as required.

(c) As $\der Gd\leq \lcs G{2^d}$, and $N_{2^d+1}=[\der Gd,G]$, we
    obtain $N_{2^d+1}\leq \der Gd\cap\lcs G{2^d+1}$. 
So it suffices to prove that $\der Gd\cap\lcs G{2^d+1}\leq N_{2^d+1}$. Assume
    now that this claim is false, and choose $x\in\der
    Gd\cap\lcs G{2^d+1}$ such that $x\not\in N_{2^d+1}$. Let $M$ be the
    subgroup $\left<N_{2^d+1},x\right>$. Note that $M$ is contained in
    $\lcs G{2^d+1}$ and that $N_{2^d+1}$ is a proper
    subgroup of $M$. If case~(i) of part~(a) is valid, then $M=\der
    Gd$. Thus, in this case, $\der Gd\leq \lcs G{2^d+1}$, and so 
Lemma~\ref{halllem} 
implies that $\der Gd/\der G{d+1}$ cannot be a
    small derived quotient. If case~(ii) of part~(a) is valid, then, as 
$\der Gd/N_{2^d+1}$ is elementary abelian with order $p^2$,  $M$ 
is a maximal subgroup in $\der Gd$, and hence $\der Gd/M$ is
    cyclic. Therefore, by Lemma~\ref{cyclemma}, $\der
    G{d+1}=[\der Gd,\der Gd]=[\der Gd,M]$. Since $M$ is contained in $\lcs
    G{2^d+1}$, $\der G{d+1}$ is contained in $[\der Gd,\lcs G{2^d+1}]$, 
which, in turn, lies in $\rc {\der Gd}G{2^d+1}=N_{2^{d+1}+1}$. 
    This, however, is a contradiction, as, 
in this case, $\der G{d+1}=N_{2^{d+1}}$. 
\end{proof}

\section{A  Lie algebra associated with $\der Gd$}\label{liealg}

We stick to the notation introduced in the previous section.

In the following two sections we assume that 
$G$ is a finite $p$-group, $\d\geq 1$,
$\der G\d/\der G{\d+1}$ is a small derived quotient, 
and that $\der G\d/[\der G\d,G]$ is elementary abelian with 
rank~2. Our aim is to show at the end of Section~\ref{3.4} that $p=2$,
and hence we prove Theorem~\ref{th1}.
We may also assume without loss of generality 
that $|\der G{\d+1}|=p$, and so $\der G{\d+1}$
lies in the centre of $G$. Lemma~\ref{2poss} shows that in this case
$|N_{2^{d+1}}|=p$,  and that $N_{2^{d+1}+1}=1$.

For the proof of Theorem~\ref{th1} we introduce three Lie algebras. The first, 
we denote it by $K$,  will be a graded Lie algebra assiciated with the 
filtration $\{N_i\}$ of $\der Gd$. The second algebra $\overline K$ will be an extension 
of $K$, in the sense that $K$ will be an ideal of $\overline K$, 
and finally, $L$ will be a small subalgebra of $\overline K$. We will 
carry out most of our computations in $L$.

Let $K$ be the graded Lie ring associated with the 
filtration $\{N_i\}$ defined in the previous section.
In other words we define $K$ as
$$
K=\bigoplus_{i=2^d}^{2^{\d+1}} N_i/N_{i+1}
$$
and the multiplication between two homogeneous elements of $K$ 
is defined as follows:
$$
[aN_{i+1},bN_{j+1}]=[a,b]N_{i+j+1}\quad\mbox{for}
\quad{a\in N_i,\ b\in N_j}.
$$
From Lemma~\ref{2poss}
it also follows that $K$ can
be viewed as a Lie algebra over $\F_p$.
Let $K_i$ denote the
homogeneous components of $K$; that is $K_i=N_i/N_{i+1}$ for $i\geq 2^\d$.
It is a consequence of Lemma~\ref{strcent} that $[K_i,K_j]\leq K_{i+j}$.

Let $Q$ denote the subgroup $[\der{G}{\d-1},G]\left(\der
G{\d-1}\right)^p$.

\begin{lemma}\label{q}
We have
\begin{enumerate}
\item[(i)] $[\der G{d-1},Q]\leq N_{2^d+1}$;
\item[(ii)] $[N_i,Q]\leq N_{i+2^{d-1}+1}$, for $i\geq 2^d$.
\end{enumerate}
\end{lemma}
\begin{proof}
(i) By~\cite[Lemma~1.1.9(ii)]{l-g}, it suffices to show that
$$
\left[\der G{d-1},\left[\der G{d-1},G\right]\right]\leq N_{2^d+1}\quad\mbox{and}
\quad \left[\der G{d-1},\left(\der
  G{d-1}\right)^p\right]\leq N_{2^d+1}.
$$
Since $[\der G{d-1},[\der G{d-1},G]]\leq\der
Gd\cap\lcs
G{2^d+1}$, using Lemma~\ref{2poss}(c), we obtain that
$[\der G{d-1},[\der G{d-1},G]]\leq
N_{2^d+1}$.
Suppose that $x,\ y\in \der G{d-1}$. Then
$[x,y^p]\equiv[x,y]^p\bmod H'^p\lcs Hp$ where $H=\left<y,[x,y]\right>$ (see equation~\eqref{coll}). 
As $\der Gd/N_{2^d+1}$ has exponent $p$, we have that $[x,y]^p\in N_{2^d+1}$,
and, since  
$$
H'^p\lcs Hp\leq \der Gd\cap \lcs G{2^d+1}=N_{2^d+1}, 
$$
we obtain
that $[x,y^p]\in N_{2^d+1}$. As $[\der G{d-1},(\der G{d-1})^p]$ is generated
as a normal subgroup by such elements $[x,y^p]$ (see~\cite[Lemma~1.1.9(iii)]{l-g}), 
we have
$[\der G{d-1},(\der G{d-1})^p]\leq N_{2^d+1}$.

(ii) As above, it suffices to show that
$$
\left[N_i,\left[\der G{d-1},G\right]\right]\leq N_{i+2^{d-1}+1}\quad\mbox{and}\quad \left[N_i,\left(\der
  G{d-1}\right)^p\right]\leq N_{i+2^{d-1}+1}.
$$
Both claims of the last displayed equation can, however, be shown as the
assertions in the proof of part~(i).
\end{proof}

Set $K_{2^{\d-1}}$ to be $\der G{\d-1}/Q$.
Clearly $K_{2^{d-1}}$ is elementary abelian, and so it can be considered as an 
$\F_p$-vector space.
Let $\overline K=K_{2^{\d-1}}\oplus K$. 
Let us define a graded Lie algebra structure on
$\overline K$. Two elements of 
$K$ are multiplied according to the product rule in $K$. If $aQ\in K_{2^{\d-1}}$
 and $b N_{i+1}\in K$ then let us define
$$
[aQ,bN_{i+1}]=[a,b]N_{2^{\d-1}+i+1}.
$$
If $aQ,\ bQ\in K_{2^{\d-1}}$
then we define 
$$
[aQ, b Q]=[a,b]N_{2^\d+1}. 
$$
We extend this product linearly and anti-symmetrically to $\overline K$.

\begin{lemma}\label{2comid}
The product on $\overline K$  is well-defined, and the algebra $\overline K$ is 
a Lie algebra over $\F_p$.
\end{lemma}
\begin{proof}
In order to show that this lemma is true we have to verify that the product
between homogeneous elements of $\overline K$ satisfies the following
properties:
\begin{enumerate}
\item[(i)] it is well-defined;
\item[(ii)] it is linear;
\item[(iii)] it is anti-symmetric;
\item[(iv)] and the Jacobi identity holds.
\end{enumerate}
Note that $K$ is a Lie ring, and so these properties obviously hold for
homogeneous elements of $K$.

(i) We claim that the product between an element of $K_{2^{d-1}}$
and a homogeneous element of $K$ is well-defined. We are required to show that
if $x,\ x'\in N_j$ and $y,\ y'\in \der G{d-1}$ such that $xN_{j+1}=x'N_{j+1}$ and
$yQ=y'Q$, then $[xN_{j+1},yQ]=[x'N_{j+1},y'Q]$. This amounts to saying that 
$$
[x,y]N_{j+2^{d-1}+1}=[x',y']N_{j+2^{d-1}+1},
$$
that is $[x,y]\equiv[x',y']\bmod N_{j+2^{d-1}+1}$. There are some $n\in
N_{j+1}$ and $g\in Q$ such that
$x'=xn$ and $y'=yg$. Now
$$
[x',y']=[xn,yg]=[x,yg][x,yg,n][n,yg]\equiv[x,yg]=[x,g][x,y][x,y,g]\equiv[x,y]
$$
where congruences are taken modulo $N_{j+2^{d-1}+1}$. Note that we used
Lemma~\ref{q}(ii). Using Lemma~\ref{q}(i), a similar argument shows that the
product is well-defined between two elements of $K_{2^{d-1}}$. 

(ii) Suppose that $x_1,\ x_2\in\der G{d-1}$, and that $y\in N_j$ with some
$j\geq 2^d$. Then, by Lemma~\ref{leftnormed}, $[x_1,y,x_2]\in
N_{j+2^{d-1}+1}$, and so
\begin{multline*}
[x_1Q+x_2Q,yN_{j+1}]=[x_1x_2Q,yN_{j+1}]=[x_1x_2,y]N_{j+2^{d-1}+1}\\
=[x_1,y][x_1,y,x_2][x_2,y]N_{j+2^{d-1}+1}=[x_1,y][x_2,y]N_{j+2^{d-1}+1}\\
=[x_1,y]N_{j+2^{d-1}+1}+[x_2,y]N_{j+2^{d-1}+1}=[x_1Q,yN_{j+1}]+[x_2Q,yN_{j+1}].
\end{multline*}
Similar argument shows that the Lie bracket is linear between homogeneous
elements of $K_{2^{d-1}}$; that
$[-xQ,yQ]=[xQ,-yQ]=-[xQ,yQ]$, for $x,\ y\in\der G{d-1}$; and that
$[-xQ,yN_{i+1}]=[xQ,-yN_{i+1}]=-[xQ,yN_{i+1}]$ for $x\in\der G{d-1}$ and $y\in
N_i$. 

(iii) It is an easy computation to show that the Lie bracket is
anti-symmetric.

(iv) Let us finally show that the Jacobi identity holds. Using the commutator
identities, we compute
\begin{equation}\label{hw1}
[x,y,z^x]=[x,y,z[z,x]]=[x,y,[z,x]][x,y,z][x,y,z,[z,x]].
\end{equation}
If $x,\ y,\ z\in\der G{d-1}$, then~\eqref{hw1} implies that
$[x,y,z^x]\equiv[x,y,z]\bmod N_{2^d+2^{d-1}+1}$. 
If $x,\ y\in\der G{d-1}$ and $z\in N_i$, then~\eqref{hw1} implies that
$[x,y,z^x]\equiv[x,y,z]\bmod N_{2^d+i+1}$. Finally, 
if $x\in \der G{d-1}$, $y\in N_i$, and $z\in N_j$, then
$[x,y,z^x]\equiv[x,y,z]\bmod N_{2^{d-1}+i+j+1}$. As the Hall-Witt identity
$$
[x,y,z^x][z,x,y^z][y,z,x^y]=1
$$
holds (see~\eqref{hw}), we find that 
$$
[x,y,z][z,x,y][y,z,x]\in N_m.
$$
where
$$
m=\left\{
\begin{array}{ll}
2^d+2^{d-1}+1 & \mbox{if $x,\ y,\ z\in \der G{d-1}$;}\\
i+2^{d}+1 & \mbox{if $x,\ y\in\der G{d-1}$ and $z\in N_i$;}\\
i+j+2^{d-1}+1 & \mbox{if $x\in\der G{d-1}, y\in N_i$, $z\in N_j$}.
\end{array}\right.
$$ 
As the Jacobi identity holds for homogeneous elements of $K$, this shows that
it also holds for homogeneous elements of $\overline K$. 
\end{proof}

Recall that
$$
\overline K=K_{2^{\d-1}}\oplus K=K_{2^{\d-1}}\oplus K_{2^\d}
\oplus K_{2^\d+1}\oplus\cdots\oplus K_{2^{\d+1}}.
$$
It is easy to see that $\overline K$ is a graded Lie algebra in the sense that
$[K_i,K_j]\in K_{i+j}$ for $i,\ j\in \{2^{\d-1},2^d,2^d+1,\ldots,2^{d+1}\}$. 
We define our object of interest: we let $L$ be the subalgebra 
$$
L=K_{2^{\d-1}}\oplus K_{2^\d}\oplus K_{2^\d+2^{\d-1}}\oplus K_{2^{\d+1}}.
$$ 
We set $L_1=K_{2^{\d-1}}$, $L_2=K_{2^\d}$, $L_3=K_{2^{\d-1}+2^\d}$ and
$L_4=K_{2^{\d+1}}$. The grading on $K$ implies that
 $L$ is also graded in the sense that
$[L_i,L_j]\leq L_{i+j}$. 
The following lemma lists the important properties of $L$.

\begin{lemma}\label{lielem}
The following is true:
\begin{enumerate}
\item[(i)] $\dim L_2 = 2$; $\dim L_3=\dim L_4=1$;
\item[(ii)] $L''=L_4$;
\item[(iii)] $L$ is generated by $L_1$.
\end{enumerate}
\end{lemma}
\begin{proof}
Claim~(i) follows because case~(ii) of part~(a) in Lemma~\ref{2poss} is
valid and because we assumed that $|N_{2^{d+1}}|=p$. 

Let us now prove~(ii) and~(iii). First we claim that 
$L_2=[L_1,L_1]$. We note that
$\der G\d=[\der G{\d-1},\der G{\d-1}]$, and so $\der G\d$ is generated by 
elements of the form $[a,b]$ with $a,\ b\in\der G{\d-1}$. Hence $L_{2^\d}$ 
is generated by elements of the form $[a,b]N_{2^\d+1}$. On the other hand,
by the definition of the product in $\overline K$, 
$[a,b]N_{2^\d+1}=[aQ,bQ]$, which implies that 
$$
L_2=K_{2^\d}=[K_{2^{\d-1}},K_{2^{\d-1}}]=[L_1,L_1].
$$

Since $\der G{\d+1}=\rc{\der G\d}G{2^d}$ and 
$|\der G{\d+1}|=p$, there are elements $a,\ b\in\der G\d$ such that
$[a,b]\neq 1$. The definition of $K$ and the fact that $\der G{\d+1}$ is
central implies
$$
[aN_{2^\d+1},bN_{2^\d+1}]=[a,b]N_{2^{\d+1}+1}=[a,b], 
$$
which, in turn, implies
that $[aN_{2^\d+1},bN_{2^\d+1}]\neq 0$. As $L_2=[L_1,L_1]$, we obtain that
$aN_{2^\d+1},\ bN_{2^\d+1}\in L'$; thus $L''\neq 0$. As $L$ has nilpotency 
class~4, and $\dim L_4=1$ we have that $L_4=L''$. This completes the proof
of~(ii). 

We have shown so far part~(ii) and that $L_1,\ L_2,\ L_4\leq\left<L_1\right>$.
It only remains to prove that $L_3\leq\left<L_1\right>$. This, however, easily
follows, because 
$$
L_4=L''=[L_2,L_2]=[[L_1,L_1],[L_1,L_1]]\leq [L_1,L_1,L_1,L_1],
$$ 
and so there must be $x_1,\
x_2,\ x_3,\ x_4\in L_1$ such that $[x_1,x_2,x_3,x_4]\neq 0$. Thus
$[x_1,x_2,x_3]\neq 0$. Since $[x_1,x_2,x_3]\in L_3$ and $\dim L_3=1$, we have
$L_3\leq \left<L_1\right>$. 
\end{proof}

\section{The proof of Theorem~\ref{th1}}\label{3.4}

We verify, in this section, that the assumptions made at the beginning of
Section~\ref{liealg} imply that $p=2$ and hence we complete the proof of
Theorem~\ref{th1}. 

We claim that elements $a,\ b,\ c\in L_1$ can 
be chosen such that the following hold:
\begin{itemize}
\item[(i)] $L_2=\left<[a,b],[a,c]\right>$ and $L_4=\left<[[a,b],[a,c]]\right>$;
\item[(ii)] $[b,c]=[a,b,a]=[a,b,b]=[a,c,a]=0$. 
\end{itemize}

By Lemma~\ref{lielem}(i), $\dim L_2=2$, and, by Lemma~\ref{lielem}(iii), $L_2=[L
_1,L_1]$, 
and so there are $x,\ y,\ u,\ v\in L_1$ such that
$L_2=\left<[x,y],[u,v]\right>$. If $[x,u]=[x,v]=[y,u]=[y,v]=0$ then
set $a=y+u$, $b=x$, $c=v$. Otherwise assume without loss of generality
that $[x,u]\neq 0$, and so there are scalars $\alpha,\ \beta\in\F_p$ such that 
$[x,u]=\alpha[x,y]+\beta[u,v]$. If $\alpha\neq 0$, then we set $a=u$, $b=x$, 
and $c=v$. If $\alpha=0$ then $\beta\neq 0$ and we set $a=x$, $b=u$, $c=y$.
It is easy to see that with each of the choices of $a,\ b,\ c$ we have that
$L_2=\left<[a,b],[a,c]\right>$ 
and that, since $L''=L_4$, it also follows that
$L_4=\left<[[a,b],[a,c]]\right>$.

We claim that $a,\ b,\ c$ above can be chosen so that $[b,c]=0$. Indeed, since
$\dim L_2=2$, some non-trivial linear combination of $[a,b]$,
$[a,c]$, $[b,c]$ is zero:
$$
\alpha[a,b]+\beta[a,c]+\gamma[b,c]=0.
$$
If $\alpha=\beta=0$ then $[b,c]=0$ and we are done. If $\alpha= 0$ and 
$\beta\neq 0$, then $[\beta a+\gamma b,c]=0$, and we replace $a$ by $\beta
  a+\gamma b$. In the new generating set $[a,c]=0$. Similarly, if $\alpha\neq
0$ and $\beta=0$, then we replace $a$ by $-\alpha a+\gamma c$.
If $\alpha\neq 0$ and $\beta\neq 0$, we replace $a$ by
  $(\beta/\alpha) a+(\gamma/\alpha) b$ and $b$ by $b+(\beta/\alpha) c$.
In the two latter cases, the new generating set satisfies
$[a,b]=0$. 
Thus after possibly reordering
the generators $a,\ b,\ c$ we may assume that $L_2=\left<[a,b],[a,c]\right>$, 
$[b,c]=0$, and $L_4=\left<[[a,b],[a,c]]\right>$. 

We claim that the elements $a,\ b,\ c$ can be chosen so that, in addition to 
what was proved in the last paragraph, we also have $[a,b,a]=[a,b,b]=
[a,c,a]=0$. First we modify the chosen elements to make $[a,b,a]$ zero. If
$[a,b,a]=0$, then we have nothing to do, and if $[a,c,a]=0$ then we 
interchange $b$ and $c$. So assume that $[a,b,a]\neq 0$ and $[a,c,a]\neq 0$.
Since $L_3$ is one-dimensional, we have that $\alpha[a,b,a]+\beta[a,c,a]=0$
for some $\alpha,\ \beta\in \F_p$ with $\alpha\beta\neq 0$. Thus
$[a,\alpha b+\beta c,a]=0$ and we replace $b$ by $\alpha b+\beta c$.

Since $[[a,b],[a,c]]=[a,b,a,c]-[a,b,c,a]=[a,b,c,a]\neq 0$ we have that 
$[a,b,c,a]\neq 0$.  Thus $a$ does not centralise $L_3$. On the other hand, 
$[a,b,a,b]=[a,b,b,a]=0$. If $[a,b,b]\neq 0$, then, as in this case
$L_3=\left<[a,b,b]\right>$, we have that $a$ centralises $L_3$, which is 
not possible. Therefore $[a,b,b]=0$ must hold. If $[a,c,a]\neq 0$ then, 
as $[a,b,c]=[a,c,b]\neq 0$ we have that 
$$
0=\alpha[a,c,a]+\beta[a,c,b]=[a,c,\alpha a+\beta b]
$$
with some $\alpha,\ \beta\in\F_p$ scalars such that $\alpha\beta\neq 0$.
If we replace $a$ by $\alpha a+\beta b$, then
we obtain the additional relation $[a,c,a]=0$. 
It is easy to see that the previously
proved relations still hold for the new elements $a,\ b,\ c$.

Hence we proved that elements $a,\ b,\ c\in L_1$ can be chosen so that
$L_2=\left<[a,b],[a,c]\right>$, $L_4=\left<[[a,b],[a,c]]\right>$, 
$[a,c]=[a,b,a]=[a,b,b]=[a,c,a]=0$. Now
\begin{multline*}
0=[a,b,a,c]+[a,c,[a,b]]+[c,[a,b],a]\\=[a,b,a,c]+[a,c,a,b]-[a,c,b,a]-[a,c,b,a]=
-2[a,c,b,a].
\end{multline*}
Thus $2[a,c,b,a]=0$, and so
$$
2[[a,b],[a,c]]=-2[a,b,c,a]=-2[a,c,b,a]=0.
$$
As $[[a,b],[a,c]]\neq 0$ we obtain that $p=2$. This completes the proof of 
Theorem~\ref{th1}.

Next we study two immediate consequences of Theorem~\ref{th1}.

\begin{lemma}\label{smalls}
Let $G$ be a finite $p$-group.
\begin{enumerate}
\item[(i)] If, for $d\geq 1$, the quotient 
$\der Gd\lcs G{2^d+1}/\lcs G{2^d+1}$ is cyclic, then, for $e\geq 1$, we have 
$\der G{d+e}\leq\lcs G{2^{d+e}+2^{e-1}}$. Further, if $\der G{d+e+1}\neq 1$,
  then 
$$
\log_p|\der G{d+e}/\der G{d+e+1}|\geq 2^{d+e}+2^{e-1}+1.
$$ 
\item[(ii)] If, for $d\geq 0$, the quotient $\der Gd\lcs G{2^d+1}/\lcs G{2^d+1}$
 can be generated by two elements, then 
$\der G{d+1}\lcs G{2^{d+1}+1}/\lcs G{2^{d+1}+1}$ is cyclic, and so, for 
$e\geq 2$, we have $\der G{d+e}\leq\lcs G{2^{d+e}+2^{e-2}}$. Further, if $\der G
{d+e+1}\neq 1$,
  then 
$$
\log_p|\der G{d+e}/\der G{d+e+1}|\geq 2^{d+e}+2^{e-2}+1.
$$ 
\end{enumerate}
\end{lemma}
\begin{proof}
(i) If $\der Gd\lcs G{2^d+1}/\lcs G{2^d+1}\cong\der Gd/(\der Gd\cap\lcs
  G{2^d+1})$ is cyclic, then, by Lemma~\ref{cyclemma},
$$
\der G{d+1}=[\der Gd,\der Gd]=[\der Gd,\der  Gd\cap\lcs G{2^d+1}]\leq
[\der Gd,\lcs G{2^d+1}]\leq\lcs G{2^{d+1}+1}.
$$
Now easy induction implies our claim.

(ii)
Suppose that $\der Gd\lcs G{2^d+1}/\lcs G{2^d+1}$ is generated by 
$x\lcs G{2^d+1}$ and $y\lcs G{2^d+1}$. If $a,\ b\in \der Gd$ then 
there are some $\alpha_a,\ \beta_a,\ \alpha_b,\ \beta_b\in\Z$ and $u_a,\
u_b\in \lcs G{2^d+1}$ such that $a=x^{\alpha_a}y^{\beta_a}u_a$ and 
$b=x^{\alpha_b}y^{\beta_b}u_b$
$$
[a,b]=[x^{\alpha_a}y^{\beta_a}u_a,x^{\alpha_b}y^{\beta_b}u_b]\equiv
[x,y]^{\alpha_a\beta_b-\beta_a\alpha_b}\bmod\lcs G{2^{d+1}+1}. 
$$
Hence 
$$
\frac{\der G{d+1}\lcs G{2^{d+1}+1}}{\lcs G{2^{d+1}+1}}\cong
\frac{\der G{d+1}}{\der G{d+1}\cap\lcs G{2^{d+1}+1}}
$$
is cyclic, and so our assertion follows from~(i). 
\end{proof}

Theorem~\ref{th1} and the previous lemma imply the following corollary.

\begin{corollary}\label{smallcor}
If $p$ is an odd prime, $G$ is a finite $p$-group, $d\geq 1$, and $\der Gd/\der 
G{d+1}$ is a small derived quotient, then, for $e\geq 1$, we have 
$\der G{d+e}\leq\lcs G{2^{d+e}+2^{e-1}}$. Further, if $\der G{d+e+1}\neq 1$, 
then
$$
\log_p|\der G{d+e}/\der G{d+e+1}|\geq 2^{d+e}+2^{e-1}+1.
$$ 
\end{corollary}
\begin{proof}
Theorem~\ref{th1} implies that $|\der Gd/[\der Gd,G]|=p$. Since $[\der
  Gd,G]\leq\der Gd\cap\lcs G{2^d+1}$, we obtain that
$$
\left|\frac{\der Gd\lcs  G{2^d+1}}{\lcs G{2^d+1}}\right|=\left|\frac{\der
  Gd}{\der Gd\cap\lcs G{2^d+1}}\right|\leq \left|\frac{\der Gd}{[\der Gd,G]}\right|=p.
$$
Hence ${\der Gd\lcs G{2^d+1}}/{\lcs G{2^d+1}}$ is cyclic and the required result follows from 
Lemma~\ref{smalls}.
\end{proof}

\section{Towards the proof of Theorem~\ref{main}}\label{secth2}

The coefficient ``3'' of the linear term of the bound in Theorem~\ref{main} is
a consequence of the following theorem.

\begin{theorem}\label{p3}
Suppose that $p$ is a prime different from $3$, 
$G$ is a finite $p$-group, $d\geq 0$, and that
$\der Gd\lcs G{2^d+1}/\lcs G{2^d+1}$ is elementary abelian. If the generator number of 
$\der G{d+2}\lcs G{2^{d+2}+1}/\lcs G{2^{d+2}+1}$ is at least~$3$, 
then either $|\lcs G{2^d}/\lcs G{2^d+1}|\geq p^4$ or $|\lcs
G{2^{d+1}+2^d}/\lcs G{2^{d+1}+2^d+1}|\geq p^2$.
\end{theorem}

The proof of this theorem  will be presented
at the end of this section.

Let $G$ denote a finite $p$-group. 
Suppose that $K$ is the subring corresponding to $\der Gd$ in the Lie 
ring associated with the lower central series of
$G$. That is
$$
K=\bigoplus_{i\geq 2^d} (\der Gd\cap\lcs Gi)\lcs G{i+1}/\lcs G{i+1}.
$$
The first homogeneous component of $K$ is 
$K_{2^d}=\der Gd\lcs G{2^d+1}/\lcs G{2^d+1}$. We  let $
L$ denote the subring
generated by $K_{2^d}$. 
For $i\geq 1$, we let $L_i=L\cap
K_{i2^d}$. Since $L$ is generated by homogeneous 
elements of the same degree, we obtain that the $L_i$ 
form a grading on $L$. 

Let me warn the reader that the Lie algebras $K$ and $L$ in the previous paragraph are different from those introduced in Section~\ref{liealg}.

\begin{lemma}\label{ls}
For $i\geq 1$, we have $L_i\cong\lcs{\der Gd}i\lcs G{i2^d+1}/
\lcs G{i2^d+1}$.
If the factor group $\der Gd\lcs G{2^d+1}/\lcs G{2^d+1}$ has exponent $p$, then $L$ can be considered as a Lie algebra over $\F_p$. Further, if
the generator number of
$\der G{d+2}\lcs G{2^{d+2}+1}/\lcs G{2^{d+2}+1}$ is at least
$3$, then the following hold:
\begin{enumerate}
\item[(i)] $\dim (L''\cap L_4)\geq 3$;
\item[(ii)] $\dim L_1\geq 3$, $\dim L_2\geq 3$, and $\dim L_3\geq 1$. 
\end{enumerate}
\end{lemma}
\begin{proof}
Since $L$ is degree-1 generated, we obtain that 
$L_i=[L_1,\ldots,L_1]$ with
$L_1$ occurring $i$ times.
Hence the definition of the product in
$L$ gives
\begin{multline*}
L_i=\big[\underbrace{L_1,\ldots,L_1}_{i\ \rm times}\big]=
\big[\underbrace{\der Gd\lcs G{2^d+1}/\lcs G{2^d+1},\ldots,\der
  Gd\lcs G{2^d+1}/\lcs G{2^d+1}}_{i\ \rm times}\big]\\=
\big[\underbrace{\der Gd,\ldots,\der Gd}_{i\ \rm times}\big]\lcs G{i2^{d}+1}/\lcs G{i2^d+1}=
\lcs{\der G{d}}i\lcs G{i2^{d}+1}/\lcs G{i2^{d}+1}. 
\end{multline*}

If $L_1=\der Gd\lcs G{2^d+1}/\lcs G{2^d+1}$ has exponent $p$, then 
$L_1$ is a vector space over $\F_p$, and, as $L$ is generated by $L_1$, we obtain that $L$ can be considered as a Lie algebra over 
$\F_p$. 

Now assume that
the generator number of 
$\der G{d+2}\lcs G{2^{d+2}+1}/\lcs G{2^{d+2}+1}$ is at least
$3$. 

(i) As $L$ is degree-1 generated, 
\begin{multline*}
L''\cap L_4=[L_2,L_2]=[[L_1,L_1],[L_1,L_1]]\\
=[[\der Gd\lcs G{2^d+1}/\lcs G{2^d+1},\der Gd\lcs G{2^d+1}/\lcs G{2^d+1}],\\{[}\der
    Gd\lcs G{2^d+1}/\lcs G{2^d+1},\der Gd\lcs G{2^d+1}/\lcs G{2^d+1}]]\\=[[\der
    Gd,\der Gd],[\der Gd,\der Gd]]\lcs G{2^{d+2}+1}/\lcs G{2^{d+2}+1}=\der
G{d+2}\lcs G{2^{d+2}+1}/\lcs G{2^{d+2}+1}.
\end{multline*}
Recall that $L_4$ is a vector space and, by assumption,  the generator number of  the factor group
$\der G{d+2}\lcs G{2^{d+2}+1}/\lcs G{2^{d+2}+1}$ is at least
$3$. Therefore $\dim (L''\cap L_4)\geq 3$.

(ii) Recall that $L$ is generated by $L_1$ and that $\dim(L''\cap L_4)=\dim [L_2,L_2]\geq 3$. 
If $L_2$ is 1-dimensional, then $[L_2,L_2]=0$, and so this is
impossible. If $L_2$ is 2-dimensional and $\{x,\ y\}$ is a basis for $L_2$
then $[L_2,L_2]$ is generated by $[x,y]$, and so $[L_2,L_2]$ is at most
1-dimensional. Hence $L_2$ is at least 3-dimensional. As $L$ is degree-1
generated, we have $L_2=[L_1,L_1]$, and so
similar argument shows
that $L_1$ is also at least 3-dimensional. Since $L_4\neq 0$, we must have
$L_3\neq 0$, and so $\dim L_3\geq 1$. 
\end{proof}

Theorem~\ref{p3} is a simple consequence of the following result.

\begin{lemma}\label{liedims}
Let $p$ be a prime different from $3$. 
If $\der Gd\lcs G{2^d+1}/\lcs G{2^d+1}$ has exponent $p$ and $\der G{d+2}\lcs G{2^{d+2}+1}/\lcs G{2^{d+2}+1}$ is a
$3$-generator abelian group, then either $\dim L_1\geq 4$ or
$\dim L_3\geq 2$. 
\end{lemma}
\begin{proof}
By Lemma~\ref{ls}, in this case, $L$ is a Lie algebra over $\F_p$, $\dim L_1\geq 3$, $\dim L_2\geq 3$, 
$\dim L_3\geq 1$, and $\dim (L_4\cap L'')\geq 3$.
Suppose that the assertion of this lemma is not valid, and so 
$\dim L_1=3$ and $\dim L_3=1$. Let 
$\{x,y,z\}$ is a basis of $L_1$. The subspace $L_4\cap
L''$ is spanned by  
$$
[[x,y],[x,z]]=[x,y,x,z]-[x,y,z,x],\quad
[[x,y],[y,z]]=[x,y,y,z]-[x,y,z,y],
$$
and
$$
[[x,z],[y,z]]=[x,z,y,z]-[x,z,z,y].
$$
If $[x,y,x]$ and $[x,y,z]$ are both zero then $[[x,y],[x,z]]=0$ and
$\dim L_4\cap L''\leq 2$. Suppose 
therefore that at least one of them is non-zero. 

Assume first that $[x,y,x]\neq 0$. Since $\dim L_3=1$ we
must have that 
$$
L_3=\bigl<[x,y,x]\bigr>\quad\mbox{and}\quad
L_4=\bigl<[x,y,x,x],[x,y,x,y],[x,y,x,z]\bigr>
$$
and there cannot be any linear dependencies among the three
generators of $L_4$. There must be some $\alpha$ such that
$[y,x,y]=\alpha[x,y,x]$. By the Jacobi identity 
$$
[y,x,[y,x]]+[x,[y,x],y]+[y,x,y,x]=0,
$$
that is $[x,y,x,y]+[y,x,y,x]=0$. In other words $[x,y,x,y]=-\alpha
[x,y,x,x]$ which means that $L_4=\bigl<[x,y,x,x],
[x,y,x,z]\bigr>$ therefore $\dim L_4\leq 2$, which is a contradiction.

Suppose now that $[x,y,z]\neq 0$. As before we must
have that 
$$
L_3=\bigl<[x,y,z]\bigr>\quad\mbox{and}\quad
L_4=\bigl<[x,y,z,x],[x,y,z,y],[x,y,z,z]\bigr>.
$$
The third homogeneous component in the free Lie algebra with rank 3 has a
basis formed by the following 8 elements:
$$
[y,x,x],\,[y,x,y],\,[y,x,z],\,[z,x,x],\,[z,x,z],\,[z,y,x],\,[z,y,y],\,[z,y,z].
$$
The relation $[z,x,y]=[y,x,z]+[z,y,x]$ also holds. Since 
the homogeneous component $L_3$ is spanned by
$-[y,x,z]=[x,y,z]$, the rest of the basis elements can be expressed as its
scalar multiples:
$$
\begin{array}{ccc}
[y,x,x]=\alpha_1[x,y,z],&[y,x,y]=\alpha_2[x,y,z],&[z,x,x]=\alpha_3[x,y,z],\\
{[z,x,z]}=\alpha_4[x,y,z],&[z,y,x]=\alpha_5[x,y,z],&[z,y,y]=\alpha_6[x,y,z],\\
&[z,y,z]=\alpha_7[x,y,z]&
\end{array}
$$
where $\alpha_i\in\F_p$ for $1\leq i\leq 7$.
The Jacobi identity implies that $[y,x,x,y]-[y,x,y,x]=0$. That
is to say that $\alpha_1[x,y,z,y]-\alpha_2[x,y,z,x]=0$. Therefore $\dim L_4\leq 
2$ unless $\alpha_1=\alpha_2=0$. In the same way,
$[z,x,x,z]-[z,x,z,x]=0$ and $[z,y,y,z]-[z,y,z,y]=0$ yield 
$\alpha_3=\alpha_4=0$ and $\alpha_6=\alpha_7=0$, respectively. It is also a
consequence of the Jacobi identity that 
$$
[y,z,[y,x]]+[z,[y,x],y]+[y,x,y,z]=0,
$$
that is 
$$
[y,z,y,x]-[y,z,x,y]-[y,x,z,y]+[y,x,y,z]=0.
$$
In other words $[y,z,x,y]+[y,x,z,y]=0$, hence
$$
\alpha_5[x,y,z,y]+[x,y,z,y]=(\alpha_5+1)[x,y,z,y]=0.
$$ 
Similarly 
$$
[y,x,[x,z]]+[x,[x,z],y]+[x,z,y,x]=0
$$
and
$$
[z,x,[z,y]]+[x,[z,y],z]+[z,y,z,x]=0
$$
yield
$$
(\alpha_5-2)[x,y,z,x]=0\quad\mbox{and}\quad (2\alpha_5-1)[x,y,z,z]=0,
$$
respectively.
Then it must follow that
$$
\alpha_5+1=0,\quad \alpha_5-2=0\quad\mbox{and}\quad 2\alpha_5-1=0
$$
which implies that $p=3$ and $\alpha_5=2$.
\end{proof}

The following example shows that the condition $p\neq 3$ is 
necessary in Lemma~\ref{liedims}.

\begin{example}
If $p=3$ then a computation with the {\sc LieNQ} program~(see Schneider~\cite{lienq}) shows
that in the Lie algebra
$$
L=\bigl<x,y,z\,|\,[y,x,x],[y,x,y],[z,x,x],[z,x,z],[z,y,x]=2[x,y,z],[z,y,y],[z,y,
z]\bigr>
$$
we have that $\dim L_1=3$, $\dim L_2=3$, $\dim L_3=1$,
$\dim L_4\cap L''=3$.
\end{example}

We can now prove Theorem~\ref{p3}.

\begin{proof}[Proof of Theorem~$\ref{p3}$]
By Lemmas~\ref{ls} and~\ref{liedims}, either
$$
\left|\frac{\lcs{\der Gd}3\lcs G{2^{d+1}+2^d+1}}{\lcs G{2^{d+1}+2^d+1}}\right|\geq p^2\quad
\mbox{or}\quad 
\left|\frac{\der Gd\lcs G{2^d+1}}{\lcs G{2^d+1}}\right|\geq p^4.
$$
Thus either 
$|\lcs G{2^{d+1}+2^d}/\lcs G{2^{d+1}+2^d+1}|\geq p^2$ or
$\quad|\lcs G{2^d}/\lcs G{2^d+1}|\geq p^4$.
\end{proof}

\begin{corollary}\label{cor}
Suppose that $p\neq 3$, $G$ is a finite $p$-group, and $d\geq 2$ 
such that the generator 
number of $\der Gi\lcs G{2^i+1}/\lcs G{2^i+1}$ is at least~$3$ for all 
$i\leq d$. Then
$\log_p|G/\der Gd|\geq 2^d+3d-2$ and $\log_p|G|\geq 2^d+3d+1$.
\end{corollary}
\begin{proof}
Note that, for $i=0,\ldots,d$, the assumptions of the theorem imply
that the order of $\lcs G{2^i}/\lcs G{2^i+1}$ is at least $p^3$. 
By Theorem~\ref{p3}, for $i=0,\ldots,d-2$, if
$\lcs G{2^i}/\lcs G{2^i+1}$ has exponent $p$, then
either $|\lcs G{2^i}/\lcs G{2^i+1}|\geq p^4$ or 
$|\lcs G{2^{i+1}+2^i}/\lcs G{2^{i+1}+2^i+1}|\geq p^2$. If the exponent of
$\lcs G{2^i}/\lcs G{2^i+1}$ is at least $p^2$, then we also have $|\lcs
G{2^i}/\lcs G{2^i+1}|\geq p^4$. Summing this up, the required result follows.
\end{proof}

\section{The proof of Theorem~\ref{main}}\label{last}

We conclude this paper with the proof of Theorem~\ref{main}.

Suppose, in this section, 
that $p\geq 5$, $G$ is a finite $p$-group, and $d$ is 
an integer such that 
$\der Gd\neq 1$. By Mann's result~\cite{Mann}, 
$\log_p|G|\geq 2^d+2d-2$, and so we may assume that $d\geq 5$. 

Let $k$ denote the smallest non-negative integer such
that the generator number of 
$\der G{k}\lcs G{2^{k}+1}/\lcs G{2^{k}+1}$
is less than 3. 
Thus the generator number of 
$\der G{i}\lcs G{2^{i}+1}/\lcs G{2^{i}+1}$ is at least~3 for all
$i\in\{0,\ldots,k-1\}$, and, in particular, $\der G{k-1}\lcs G{2^{k-1}+1}/\lcs
G{2^{k-1}+1}$ has order at least $p^3$.
This also implies that 
$\der G{k-1}\neq 1$, and so we may assume without loss of generality that $d\geq k-1$. 
If $k=0$, then $G$ is generated by two elements, and so, 
$\log_p|\der Gi/\der G{i+1}|\geq 2^i+2^{i-2}+1$ for $i\in\{2,\ldots,d-1\}$. As
$\der Gd\neq 1$, we have $\log_p|G|\geq 2^d+2^{d-2}+d-1$. As $d\geq 5$, we
obtain that the assertion of the theorem is true in this case. If $k=1$, then,
for $i=3,\ldots,d-1$,  
we have  $\log_p|\der Gi/\der G{i+1}|\geq 2^i+2^{i-3}+1$. Then $\log_p|G|\geq 
2^d+2^{d-3}+d$, and so the theorem is valid also in this case. 
If $k=2$, then $\log_p|\der Gi/\der G{i+1}|\geq 2^i+2$ for $i=0,\ 1$,
$|\der G2/\der G3|\geq p^5$, $|\der G3/\der G4|\geq p^9$. However, 
by Corollary~\ref{smallcor} 
we have that either $|\der G2/\der G3|\geq p^6$ or 
$|\der G3/\der G4|\geq p^{10}$. For $i\geq 4$ we have that 
$\log_p|\der Gi/\der G{i+1}|\geq 2^i+2^{i-4}+1$. 
Thus $\log_p|G|\geq 2^d+2^{d-4}+d+2$. Therefore, in this case, we have 
nothing to prove.
Hence we may assume without loss of generality that $k\geq 3$. 

By Corollary~\ref{cor}, 
$$
\log_p|G/\der G{k-1}|\geq 2^{k-1}+3k-5\quad\mbox{and}\quad\log_p|G|\geq 
2^{k-1}+3k-2. 
$$
Thus, if $d=k-1$, then $\log_p|G|\geq 2^d+3d+1$, the result is valid.

Let $d\geq k$. It follows   from Lemma~\ref{2poss}(c) that if $\log_p|\der
G{k-1}/\der G{k}|=2^{k-1}+1$, then $[\der G{k-1},G]=\der G{k-1}\cap \lcs
G{2^{k-1}+1}$, and so
$$
\frac{\der G{k-1}\lcs G{2^{k-1}+1}}{\lcs G{2^{k-1}+1}}\cong\frac{\der G{k-1}}
{\der G{k-1}\cap\lcs G{2^{k-1}+1}}=\frac{\der G{k-1}}{[\der G{k-1},G]}.
$$
However, in this case, $|\der G{k-1}/[\der G{k-1},G]|\leq p^2$
(see Lemma~\ref{2poss}), and so we must have that $\log_p|\der G{k-1}/\der
G{k}|\geq 2^{k-1}+2$.
Thus
\begin{multline*}
\log_p|G/\der G{k}|=\log_p|G/\der G{k-1}|+\log_p|\der G{k-1}/\der G{k}|\\
\geq 2^{k-1}+3k-5+2^{k-1}+2=2^{k}+3k-3.
\end{multline*}
Thus, if $d=k$, then $|\der G{k}|\geq p$, and we find
$$
\log_p|G|=2^{k}+3k-2=2^d+3d-2,
$$
and the result is, again, valid.

Assume now that $d\geq k+1$. As $\log_p|\der G{k}/\der G{k+1}|\geq
2^{k}+1$, we have
\begin{multline*}
\log_p|G/\der G{k+1}|=\log_p|G/\der G{k}|+\log_p|\der G{k}/\der G{k+1}|
\\\geq 2^{k}+3k-3+2^{k}+1=2^{k+1}+3k-2.
\end{multline*}
Therefore, if $d=k+1$, then $\log_p|G|\geq 2^d+3d-4$, and so the theorem is
valid also in this case. 

Therefore we suppose that $d\geq k+2$. 
In this case, $\log_p|\der G{k+1}/\der
G{k+2}|\geq 2^{k+1}+1$. However, as $k\geq 3$, 
by Corollary~\ref{smallcor}, either 
$$
\log_p|\der G{k}/\der G{k+1}|\geq
2^{k}+2 
$$
or 
$$
\log_p|\der G{k+1}/\der
G{k+2}|\geq 2^{k+1}+2
$$
Thus 
\begin{multline*}
\log_p|G/\der G{k+2}|=\log_p|G/\der G{k+1}|+\log_p|\der G{k+1}/\der G{k+2}|
\\\geq 2^{k+1}+3k-2+2^{k+1}+2=2^{k+2}+3k.
\end{multline*}
Hence if $d=k+2$, 
then $\log_p|G|\geq 2^{k+2}+3k+1=2^d+3d-5$. Therefore, in
this case, the theorem is valid.

Assume now that $d\geq k+3$. Corollary~\ref{smallcor} yields 
$\log_p|\der G{k+2}/\der G{k+3}|\geq 2^{k+2}+2$. Hence 
\begin{multline*}
\log_p|G/\der G{k+3}|=\log_p|G/\der G{k+2}|+\log_p|\der G{k+2}/\der G{k+3}|
\\\geq 2^{k+2}+3k+2^{k+2}+2=2^{k+3}+3k+2.
\end{multline*}
Hence if $d=k+3$, then $\log_p|G|\geq 2^{k+3}+3k+3=2^d+3d-6$, and so the
result holds. 

If $d\geq k+4$, then, by Lemma~\ref{smallcor}, $\log_p|\der Gi/\der
G{i+1}|\geq 2^i+3$, for all $i\in\{k+3,\ldots,d-1\}$. 
Thus in this case the theorem also holds.

\section*{Acknowledgment}
Much of the research presented in this paper was carried out while I was a PhD
student at The Australian National University. I am particularly grateful to
my PhD supervisor, Mike Newman, for his continuous support.

\bibliographystyle{alpha}

\end{document}